\documentclass[a4paper]{amsart}
\usepackage{graphicx}
\usepackage{amsmath}
\usepackage{latexsym}
\usepackage{amssymb}
\usepackage[all]{xy}
\xyoption{matrix}
\xyoption{arrow}
\begin{document}

\renewcommand{\th}{\operatorname{th}\nolimits}
\newcommand{\rej}{\operatorname{rej}\nolimits}
\newcommand{\extto}{\xrightarrow}
\renewcommand{\mod}{\operatorname{mod}\nolimits}
\newcommand{\ul}{\underline}
\newcommand{\Sub}{\operatorname{Sub}\nolimits}
\newcommand{\ind}{\operatorname{ind}\nolimits}
\newcommand{\Fac}{\operatorname{Fac}\nolimits}
\newcommand{\add}{\operatorname{add}\nolimits}
\newcommand{\soc}{\operatorname{soc}\nolimits}
\newcommand{\Hom}{\operatorname{Hom}\nolimits}
\newcommand{\Rad}{\operatorname{Rad}\nolimits}
\newcommand{\RHom}{\operatorname{RHom}\nolimits}
\newcommand{\uHom}{\operatorname{\underline{Hom}}\nolimits}
\newcommand{\End}{\operatorname{End}\nolimits}
\renewcommand{\Im}{\operatorname{Im}\nolimits}
\newcommand{\Ker}{\operatorname{Ker}\nolimits}
\newcommand{\Coker}{\operatorname{Coker}\nolimits}
\newcommand{\Ext}{\operatorname{Ext}\nolimits}
\newcommand{\op}{{\operatorname{op}}}
\newcommand{\Ab}{\operatorname{Ab}\nolimits}
\newcommand{\id}{\operatorname{id}\nolimits}
\newcommand{\pd}{\operatorname{pd}\nolimits}
\newcommand{\A}{\operatorname{\mathcal A}\nolimits}
\newcommand{\C}{\operatorname{\mathcal C}\nolimits}
\newcommand{\D}{\operatorname{\mathcal D}\nolimits}
\newcommand{\X}{\operatorname{\mathcal X}\nolimits}
\newcommand{\Y}{\operatorname{\mathcal Y}\nolimits}
\newcommand{\F}{\operatorname{\mathcal F}\nolimits}
\newcommand{\Z}{\operatorname{\mathbb Z}\nolimits}
\renewcommand{\P}{\operatorname{\mathcal P}\nolimits}
\newcommand{\T}{\operatorname{\mathcal T}\nolimits}
\newcommand{\G}{\Gamma}
\renewcommand{\L}{\Lambda}
\newcommand{\bdot}{\scriptscriptstyle\bullet}
\renewcommand{\r}{\operatorname{\underline{r}}\nolimits}
\newtheorem{lemma}{Lemma}[section]
\newtheorem{prop}[lemma]{Proposition}
\newtheorem{cor}[lemma]{Corollary}
\newtheorem{theorem}[lemma]{Theorem}
\newtheorem{remark}[lemma]{Remark}
\newtheorem{definition}[lemma]{Definition}
\newtheorem{example}[lemma]{Example}

\newcommand{\im}{\mbox{im}}

\newtheorem{ex}[lemma]{Example}
\newtheorem{conj}{Conjecture}

\title[A geometric description of the $m$-cluster categories of type $D_n$]{A geometric description of the $m$-cluster categories of type $D_n$}

\author[Baur]{Karin Baur}
\address{Department of Mathematics \\
University of Leicester \\
University Road \\
Leicester LE1 7RH \\
England
}
\email{k.baur@mcs.le.ac.uk}

\author[Marsh]{Robert J. Marsh}
\address{Department of Pure Mathematics \\
University of Leeds \\
Leeds LS2 9JT \\
England
}
\email{marsh@maths.leeds.ac.uk}


\begin{abstract}
We show that the $m$-cluster category of type $D_n$ is equivalent to a
certain geometrically-defined category of arcs in a punctured
regular $nm-m+1$-gon. This generalises a result of Schiffler for $m=1$.
We use the notion of the $m$th power of a translation quiver to realise
the $m$-cluster category in terms of the cluster category.
\end{abstract}

\maketitle

\section*{Introduction}

Let $k$ be a field and $Q$ a quiver of Dynkin type $\Delta$.
Let $D^b(kQ)$ denote the
bounded derived category of finite dimensional $kQ$-modules. Let $\tau$
denote the Auslander-Reiten translate of $D^b(kQ)$ and let $S$ denote
the shift. For $m\in \mathbb{N}$ the \emph{$m$-cluster category} associated
to $kQ$ is the orbit category
$$\mathcal{C}^m_{\Delta}:=\frac{D^b(kQ)}{S^m\tau^{-1}}.$$
This category was introduced in~\cite{keller} and has been studied
by the authors~\cite{baurmarsh}, Thomas~\cite{thomas},
Wralsen~\cite{wralsen} and Zhu~\cite{zhu}.
It is known that $\mathcal{C}^m_{\Delta}$ is triangulated~\cite{keller},
Krull-Schmidt and has almost split triangles~\cite[1.2,1.3]{bmrrt}.

The $m$-cluster category is a generalisation of the cluster category.
The cluster category was introduced in~\cite{ccs1} (for type $A$)
and~\cite{bmrrt}
(general hereditary case), and can be regarded as the case $m=1$ of
the $m$-cluster category.
Keller has shown that the $m$-cluster category is Calabi-Yau
of dimension $m+1$~\cite{keller}. We remark that such Calabi-Yau categories
have also been studied in~\cite{kellerreiten}. One of the aims of the
definition of the cluster category was to model the Fomin-Zelevinsky
cluster algebra~\cite{fominzelevinsky} representation-theoretically.

We show that $\mathcal{C}^m_{D_n}$
can be realised geometrically in terms of a category of arcs in a
punctured polygon with $nm-m+1$ vertices. This generalises a result of
Schiffler~\cite{schiffler}, who considered the case $m=1$. We remark
that the punctured polygon model for the cluster
algebra of type $D_n$
appears in work of Fomin, Schapiro and Thurston~\cite{fst} as part of a
more general set-up, building on~\cite{fg1,fg2,gsv1,gsv2} which consider
links between cluster algebras and Teichm\"{u}ller theory.

Also, such a geometric realisation of a cluster
category first appeared (with a construction for type $A_n$ in the
case $m=1$) in~\cite{ccs1}.

Our approach is based on the idea of the \emph{$m$th power} of a translation
quiver introduced in~\cite{baurmarsh}. We show that, with a slight modification
of the definition for $m=2$, the Auslander-Reiten quiver of
$\mathcal{C}^m_{D_n}$ can be realised as a connected component of the $m$th
power of the Auslander-Reiten quiver of 
$\mathcal{C}^1_{D_{nm-m+1}}$. 
In Section~\ref{se:toralexample} we show that, if this 
modification is not
made, the square of the Auslander-Reiten quiver of $\mathcal{C}^1_{D_4}$
has a connected component whose underlying topological space is a torus.

\section{Notation and Definitions} \label{notation}

\begin{figure}[ht] 
\begin{center}
\includegraphics{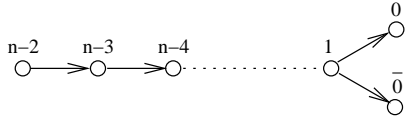} 
\caption{Quiver of type $D_n$}\label{fig1}
\label{dnquiver}
\end{center}
\end{figure}

Let $Q$ be a quiver of underlying Dynkin type $D_n$. 
The vertices 
of $Q$ are labelled $0,\overline{0},1,\dots,n-2$ and the 
arrows are $i\to i-1$ ($i=1,\dots,n-2$) together with 
$1\to \overline{0}$; see Figure~\ref{dnquiver}.

\vspace{.5cm}
We now recall the Auslander-Reiten quiver of the 
cluster category $\mathcal{C}_{D_n}$ (see~\cite[\S1]{bmrrt},~\cite{happel}).
It is a stable translation quiver built from $n$ copies 
of $Q$. We denote it by 
$\Gamma(D_n,1)$. 
The vertices of $\Gamma(D_n,1)$ are 
$V(D_n,1):=\mathbb{Z}_n\times
\{0,\overline{0},1,\dots,n-2\}$. 
The arrows are 
\[
\left.
\begin{array}{l}
(i,j)\to (i,k) \\ 
(i,k)\to (i+1,k)
\end{array}
\right\}
\quad \text{whenever there is an arrow $j\to k$ in $D_{n}$.}
\]
Finally, the translation $\tau$ is given by 
\[
\tau(i,j)=\left\{\begin{array}{ll} 
(i-1,\overline{j}), & \text{if $i=0,j\in\{0,\overline{0}\}$ and 
$n$ is odd,}\\
(i-1,j), & \text{otherwise.}
\end{array}\right.
\]
We use the convention that $\overline{\overline{0}}=0$. 
Note that the switch described here only occurs 
for odd $n$.

\begin{figure}
\begin{center}
\includegraphics{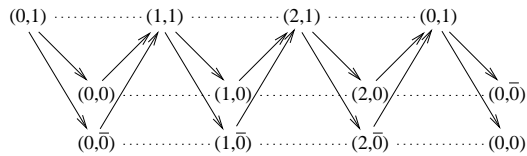} 
\caption{The quiver $\Gamma(D_3,1)$}
\label{fi:quiverd3}
\end{center}
\end{figure}

\begin{figure}
\begin{center}
\includegraphics{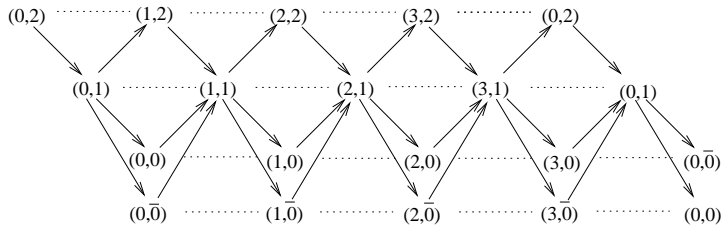} 
\caption{The quiver $\Gamma(D_4,1)$}
\label{fi:quiverd4}
\end{center}
\end{figure}

As an example, we draw the 
quivers $\Gamma(D_n,1)$ for 
$n=3$ and $n=4$; see Figures~\ref{fi:quiverd3} and~\ref{fi:quiverd4}.
The translation $\tau$ is indicated by dotted lines
(it is directed to the left). 

We recall the notion of the $m$-th power of a translation quiver 
(cf.~\cite{baurmarsh}). If $\Gamma$ is a translation 
quiver with translation $\tau$, then the 
quiver $\Gamma^m$, the {\em $m$-th power of $\Gamma$}, is the 
quiver whose objects are 
the same as the objects of $\Gamma$ and whose arrows 
are the sectional paths (in $\Gamma$) of length $m$. 
A path $x=x_0\to x_1\to\dots\to x_{m-1}\to x_m=y$ is 
said to be {\em sectional} if $\tau x_{i+1}\neq x_{i-1}$ for 
$i=1,\dots,m-1$ (in the cases where $\tau x_{i+1}$ is defined), 
cf.~\cite{ringel}. 

One of our goals is to realise the Auslander-Reiten quiver for
the $m$-cluster category of type $D_n$ in terms of the $m$th power
of the Auslander-Reiten quiver of a cluster category of type 
$D_{nm-m+1}$.
To be able to do this, we introduce a new class of sectional paths. 
\begin{definition} \rm
Let $\Gamma=\Gamma(D_n,1)$ be the translation quiver defined above,
with vertices
$V(D_n,1)=\mathbb{Z}_n\times\{0,\overline{0},1,\dots,n-2\}$. 
We say that 
a sectional path $x=x_0\to x_1\to\dots\to x_{m-1}\to x_m=y$ (where 
$x_i\in V(D_n,1)$) is {\em restricted} if there is no $i$ such that
$x_{i+1}=(r,0)$ for some $r$ while $x_{i-1}=(r-1,\overline{0})$ 
or such that $x_{i+1}=(r,\overline{0})$ for some $r$ 
while $x_{i-1}=(r-1,0)$. 
\end{definition}

\begin{remark} \rm \label{re:restricted}
Note that unless $m=2$, the restricted sectional paths 
of length $m$ in $\Gamma$ are exactly the sectional paths 
of length $m$. We can see this as follows. Firstly, it is clear
that any sectional path of length $1$ is necessarily restricted.
Suppose that $m>2$. Let 
$x=x_0\to x_1\to\dots\to x_{m-1}\to x_m=y$ be sectional,
and suppose that there is an $i\in\{1,\dots,m-1\}$ such 
that $x_{i+1}=(r,0)$ and $x_{i-1}=(r-1,\overline{0})$. 
Then $x_i=(r,1)$. In case $i=1$ we have $x_{i+2}=(r+1,1)$, and
it follows that the original path is not sectional, a contradiction.
Similarly, if $i>1$, we have $x_{i-2}=(r-1,1)$, and again the
original path is not sectional.

Hence the only sectional paths that are not restricted 
are the paths of the form 
$(i,0)\to (i,1) \to (i+1,\overline{0})$ and 
$(i,\overline{0})\to (i,1)\to (i+1,0)$ ($i\in\mathbb{Z}_n$). 
\end{remark}

With this new notion we are now ready to introduce a 
restricted version of the translation quiver $((\Gamma(D_n,1))^m,\tau^m)$. 
We define a translation quiver \linebreak
$(\mu_m(\Gamma(D_n,1)),\tau^m)$ as follows. 
The vertices of $(\mu_m(\Gamma(D_n,1)),\tau^m)$ are the same as the vertices 
of $(\Gamma(D_n,1),\tau^m)$, i.e. 
$\mathbb{Z}_n\times\{0,\overline{0},1,\dots,n-2\}$, 
the arrows are the restricted sectional paths of length $m$ 
in $(\Gamma(D_n,1),\tau^m)$ and the translation is 
$\tau^m$. 

\begin{lemma}
For any $m$, the pair $(\mu_m(\Gamma(D_n,1)),\tau^m)$ is a 
stable translation quiver.
\end{lemma}
\begin{proof}
We firstly note that the unrestricted version, 
$(\Gamma(D_n,1)^m,\tau^m)$, is a stable translation 
quiver by~\cite[6.2]{baurmarsh}. 
By Remark~\ref{re:restricted}, the quiver 
$(\mu_m(\Gamma(D_n,1)),\tau^m)$
is the same as $(\Gamma(D_n,1)^m,\tau^m)$ if $m\not=2$, so the result
follows in this case.

Now assume that $m=2$ and fix a vertex $x$ in $\Gamma(D_n,1)$. 
To show that \linebreak
$(\mu_m(\Gamma(D_n,1)),\tau^m)$ is a translation quiver, we need
to show that there is a restricted sectional path of length $2$ from $y$ to
$x$ if and only if there is a restricted sectional path of length $2$ from
$\tau^2(x)$ to $y$. Since the restricted sectional paths in $\Gamma(D_n,1)$ of
length $2$ starting or ending at $x$ are the same as the sectional paths
provided $x$ is not of the form $(i,0)$ or $(i,\overline{0})$, we are
reduced to this case. If $x=(i,0)$, the sectional paths of length $2$
ending in $x$ are $(i,2)\rightarrow (i,1)\rightarrow (i,0)$ and
$(i-1,\overline{0})\rightarrow (i,1)\rightarrow (i,0)$. The sectional
paths of length $2$ starting at $\tau^2(x)=(i-2,0)$ are
$(i-2,0)\rightarrow (i-1,1)\rightarrow (i,2)$ and
$(i-2,0)\rightarrow (i-1,1)\rightarrow (i-1,\overline{0})$.
The second path only in each case is not restricted, so we see that there is a
restricted sectional path of
length $2$ from $y$ to $x$ if and only if there is a restricted sectional
path of length $2$ from $\tau^2(x)$ to $y$. The argument in case
$x=(i,\overline{0})$ is similar. Hence $(\mu_m(\Gamma(D_n,1)),\tau^m)$ is a 
translation quiver. 

By construction, no vertex is projective 
and $\tau^m$ is defined on all vertices (since $\tau$ 
is). Therefore, $(\mu_m(\Gamma(D_n,1)),\tau^m)$ is stable.
\end{proof}
%
%

%
%

\section{The $m$-cluster category of type $D_n$ as a component of a restricted
$m$th power}
\label{se:mthpower}

Let $n,m\in\mathbb{N}$, with $n\geq 3$.
We recall that~\cite{happel} the derived category of a quiver of
Dynkin type $D_n$ has vertices $\mathbb{Z}\times \{0,\overline{0},1,2,\ldots
,n-2\}$ and arrows given by
$(i,j)\rightarrow (i,j-1)$ and
$(i,j-1)\rightarrow (i+1,j)$ for $1\leq j\leq n-2$,
and $(i,1)\rightarrow (i,\overline{0})$ and
$(i,\overline{0})\rightarrow (i+1,1)$,
where $i\in\mathbb{Z}$ is arbitrary.
We also have that
$$S^m(i,0)=\left\{ \begin{array}{cc} (i+m,0), & \text{$nm$ even,} \\
(i+m,\overline{0}), & \text{$nm$ odd.} \end{array}\right., $$ 
while $S^m(i,j)=(i+m,j)$, otherwise.

Let $\Gamma(D_n,m)$ be the quiver with vertices
$$V(D_n,m)=\{(i,j):i\in\mathbb{Z}_{nm-m+1},\ 
j\in \{0,\overline{0},1,2,\ldots, n-2\}\}.$$
The arrows are given by $(i,j)\rightarrow (i,j-1)$ and
$(i,j-1)\rightarrow (i+1,j)$ for $1\leq j\leq n-2$,
and $(i,1)\rightarrow (i,\overline{0})$ and
$(i,\overline{0})\rightarrow (i+1,1)$,
where $i\in\mathbb{Z}_{nm-m+1}$ is arbitrary and the addition is modulo
$nm-m+1$.
We also define
$$\widetilde{\tau}(i,j)=
\left\{ \begin{array}{cc} (i-1,\overline{j}) & \mbox{\ if\ }
i=0,\ j\in\{0,\overline{0}\}\mbox{\ and\ }nm\mbox{\ is\ odd,} \\
(i-1,j) & \mbox{otherwise.}
\end{array}\right.$$

It follows from the construction of $\mathcal{C}^m_{D_n}$
and the above description of the derived category that
$(\Gamma(D_n,m),\widetilde{\tau})$ is the
Auslander-Reiten quiver of $\mathcal{C}^m_{D_n}$
(and, in particular, is a stable translation quiver).

The vertices of the Auslander-Reiten quiver $\Gamma(D_{nm-m+1},1)$
of $\mathcal{C}^1_{D_{nm-m+1}}$ are
$$V(D_{nm-m+1},1)=
\mathbb{Z}_{nm-m+1}\times \{0,\overline{0},1,2,\ldots ,nm-m-1\}.$$
The arrows are given by $(i,j)\rightarrow (i,j-1)$ and
$(i,j-1)\rightarrow (i+1,j)$ for $1\leq j\leq nm-m-1$,
and $(i,1)\rightarrow (i,\overline{0})$ and
$(i,\overline{0})\rightarrow (i+1,1)$,
where $i$ is arbitrary and the addition is modulo $nm-m+1$.
We also have
$$\tau(i,j)=\left\{ \begin{array}{cc} (i-1,\overline{j}) & \mbox{\ if\ }
i=0,\ j\in\{0,\overline{0}\}\mbox{\ and\ }nm-m+1\mbox{\ is\ odd,} \\
(i-1,j) & \mbox{otherwise.}
\end{array}\right.$$

\begin{definition}
We define a map $\sigma'$ from $V(D_n,m)$ to $V(D_{nm-m+1},1)$ as follows.
We set $\sigma'(i,j)=(im,jm)$ whenever $j\not\in\{0,\overline{0}\}$ or
$j\in\{0,\overline{0}\}$, $m$ is odd and $n$ is even. Otherwise, we
have $j=0$ or $\overline{0}$ and we set
$$\sigma'(i,j)=\left\{\begin{array}{cc}
(im,jm), & \left\lfloor\frac{im}{nm-m+1}\right\rfloor \mbox{\ even,} \\
(im,\overline{jm}), & \left\lfloor\frac{im}{nm-m+1}\right\rfloor \mbox{\ odd.}
\end{array}\right.$$
Here we restrict $i$ to lie in the set $\{0,1,2, n-1\}$.
\end{definition}

We use the usual convention that, for a real number $x$, $\lfloor x\rfloor$
denotes the largest integer $k$ such that $k\leq x$.

Let $V:=\{(r,s)\in V(D_{nm-m+1},1)\,:\,m|s\}$. Here we adopt the convention
that $m|\overline{0}$.

\begin{lemma}
With $\sigma'$ defined as above, we have that $\im(\sigma')=V$.
\end{lemma}

\begin{proof}
First we note that it is clear from the definition of $\sigma'$ that
$\im(\sigma')\subseteq V$. Let $(r,s)\in V(D_{nm-m+1},1)$ and suppose that
$m|s$. Suppose first that $s\not=0,\overline{0}$.
Write $s=km$ for $k\in\mathbb{Z}$. We have
$r=r(nm-m+1-(n-1)m)=-r(n-1)m$ in $\mathbb{Z}_{nm-m+1}$, and it follows that 
$\sigma'(-r(n-1),k)=(-r(n-1)m,km)=(r,s)$ so
$(r,s)\in\im(\sigma')$. If $s=0$ or $\overline{0}$ then
$\{\sigma'(-r(n-1),0),\sigma'(-r(n-1),\overline{0})\}=
\{(r,s),(r,\overline{s})\}$ and we are done.
\end{proof}

Let $\Gamma$ denote the full subquiver of $\mu_m(\Gamma(D_{nm-m+1},1))$
induced by $V$, and let $\sigma$ be the (surjective) map obtained by
restricting the codomain of $\sigma'$ to $V$. We will show that $\sigma$
is an isomorphism from $\Gamma(D_n,m)$ to $\Gamma$ and that
$\Gamma$ is a connected component of $\mu_m(\Gamma(D_{nm-m+1},1))$.

\begin{lemma} \label{closed}
Let $x:=(r,s)\in V$ and suppose that
$$x=x_0\rightarrow x_1\rightarrow x_2\rightarrow \cdots \rightarrow
x_m=y$$
is a restricted sectional path of length $m$ in $\Gamma(D_{nm-m+1},1)$.
Then $y\in V$.
\end{lemma}

\begin{proof}
If $s\not=0,\overline{0}$ we can argue as in~\cite[7.1]{baurmarsh} to
see that $y$ is either $(r,s-m)$, $(r+m,s+m)$ or $(r,\overline{s-m})$,
where the last possibility only arises if $s=m$.
If $s=0$ or $\overline{0}$ a similar argument shows that $y$ must
be $(r+m,m)$. If $m=2$ we have the sectional paths
$(r,0)\rightarrow (r+1,1)\rightarrow (r+2,\overline{0})$ and
$(r,\overline{0})\rightarrow (r+1,1)\rightarrow (r+2,0)$, but, by definition,
these are not restricted sectional paths.
\end{proof}

\begin{lemma}
The map $\sigma:\Gamma(D_n,m)\rightarrow \Gamma$ is an isomorphism of quivers.
\end{lemma}

\begin{proof}
Since $|V|=|V(D_n,m)|$ and $\sigma$ is surjective, it
follows that $\sigma$ is bijective.
The arrows in $\Gamma(D_n,m)$ are of the form
$(i,j)\rightarrow (i,j-1)$, $(i,j)\rightarrow (i+1,j+1)$,
$(i,1)\rightarrow (i,\overline{0})$ or $(i,\overline{0})\rightarrow (i+1,1)$.
The arrows in $V\subseteq \mu_m(\Gamma(D_{nm-m+1},1))$ are of the form
$(r,s)\rightarrow (r,s-m)$, $(r,s)\rightarrow (r+m,s+m)$,
$(r,m)\rightarrow (r,\overline{0})$ or $(r,\overline{0})\rightarrow (r+m,m)$
(see the proof of Lemma~\ref{closed}).
It follows that $\sigma$ is an isomorphism of quivers.
\end{proof}

\begin{prop}
The map $\sigma:\Gamma(D_n,m)\rightarrow \Gamma$ is an isomorphism
of translation quivers. Its image, $\Gamma$, is a connected component of
$\mu_m(\Gamma(D_{nm-m+1},1))$.
\end{prop}

\begin{proof}
By Lemma~\ref{closed}, both statements of the theorem will follow if we can
show that, for all $(i,j)\in V(D_n,m)$, $\sigma(\widetilde{\tau}(i,j))=
\tau^m(\sigma(i,j))$, since this will also imply that the image
of $\sigma$ is closed under $\tau^m$.
We firstly note that if $j\not=0,\overline{0}$ then
$\widetilde{\tau}(i,j)=(i-1,j)$ while $\tau^m(im,jm)=((i-1)m,jm)$.
Since $\sigma(i,j)=(im,jm)$ and $\sigma(i-1,j)=((i-1)m,jm)$,
the result holds. So we are left with the case where $j=0$ or $\overline{0}$.
We break this down into cases, considering first the case where $j=0$.

\noindent {\bf Case (a)}: $m$ odd and $n$ even. \\
In this case we have that $nm-m+1$ is even, and $nm$ is even, so for any
$(i,0)\in V(D_n,m)$, $\tau^m(im,0)=((i-1)m,0)$ while
$\widetilde{\tau}(i,0)=(i-1,0)$. Since $\sigma(i-1,0)=((i-1)m,0)$
and $\sigma(i,0)=(im,0)$, we are done.

\noindent {\bf Case (b)}: $m$ is even. \\
In this case we have that $nm-m+1$ is odd, so for $l=0$ or $\overline{0}$,
we have:
$$\tau^m(im,l)=\left\{ \begin{array}{cc} ((i-1)m,\overline{l}), & im\mod
nm-m+1 \in\{0,1,\ldots ,m-1\}, \\
((i-1)m,l), & \mbox{otherwise.}
\end{array}\right.
$$
Since $m$ is even, $nm$ is even, so
$\widetilde{\tau}(i,0)=(i-1,0)$ for all $i$.

(i) Suppose first that $im\not\in \{0,1,\ldots ,m-1\}$. Then
$\left\lfloor\frac{im}{nm-m+1}\right\rfloor=\left\lfloor\frac{(i-1)m}{nm-m+1}\right\rfloor$.
It follows that either $\sigma(i-1,0)=((i-1)m,0)$ and $\sigma(i,0)=(im,0)$
or $\sigma(i-1,0)=((i-1)m,\overline{0})$ and $\sigma(i,0)=(im,\overline{0})$.
In either case we see that $\sigma(\widetilde{\tau}(i,0))=\tau^m(\sigma(i,0))$.

(ii) Suppose next that $im\in\{1,\ldots ,m-1\}$. Then
$\left\lfloor\frac{im}{nm-m+1}\right\rfloor-1=\left\lfloor\frac{(i-1)m}{nm-m+1}\right\rfloor$.
It follows that either $\sigma(i-1,0)=((i-1)m,\overline{0})$ and
$\sigma(i,0)=(im,0)$
or $\sigma(i-1,0)=((i-1)m,0)$ and $\sigma(i,0)=(im,\overline{0})$.
In either case we see that $\sigma(\widetilde{\tau}(i,0))=\tau^m(\sigma(i,0))$.

(iii) Finally, suppose that $i=0$. Then $i-1\equiv (n-1)m\mod nm-m+1$, and
$\left\lfloor\frac{im}{nm-m+1}\right\rfloor=0$ is even while
\begin{eqnarray*}
\left\lfloor\frac{(n-1)m}{nm-m+1}\right\rfloor & = &
\left\lfloor\frac{(n-1)m^2}{nm-m+1}\right\rfloor \\
& = & \left\lfloor\frac{(nm-m+1)m}{nm-m+1}-\frac{m}{nm-m+1}\right\rfloor \\
& = & \left\lfloor m-\frac{m}{nm-m+1}\right\rfloor=m-1,
\end{eqnarray*}
is odd (using here the fact that $m<(n-1)m+1=nm-m+1$).
It follows that $\sigma(i-1,\overline{0})=((i-1)m,0)$ and
$\sigma(i,0)=(im,0)$ and thus that
$\sigma(\widetilde{\tau}(i,0))=\tau^m(\sigma(i,0))$.

\noindent {\bf Case (c)}: $n,m$ both odd. \\
In this case we have that $nm-m+1$ is odd, so for $l=0$ or $\overline{0}$,
we have:
$$\tau^m(im,l)=\left\{ \begin{array}{cc} ((i-1)m,\overline{l}) & im\in\{0,1,\ldots ,m-1\}, \\
((i-1)m,l) & \mbox{otherwise.}
\end{array}\right.
$$
Since $n$ and $m$ are both odd, $nm$ is odd, so
$$\widetilde{\tau}(i,0)=
\left\{ \begin{array}{cc} ((i-1),\overline{0}) & i=0, \\
(i-1,0) & \mbox{otherwise.}
\end{array}\right.
$$
(i) Suppose first that $im\not\in \{0,1,\ldots ,m-1\}$. Then
$\left\lfloor\frac{im}{nm-m+1}\right\rfloor=\left\lfloor\frac{(i-1)m}{nm-m+1}\right\rfloor$.
It follows that either $\sigma(i-1,0)=((i-1)m,0)$ and $\sigma(i,0)=(im,0)$
or $\sigma(i-1,0)=((i-1)m,\overline{0})$ and $\sigma(i,0)=(im,\overline{0})$.
In either case we see that $\sigma(\widetilde{\tau}(i,0))=\tau^m(\sigma(i,0))$.

\noindent (ii) Suppose next that $im\in\{1,\ldots ,m-1\}$. Then
$\left\lfloor\frac{im}{nm-m+1}\right\rfloor-1=\left\lfloor\frac{(i-1)m}{nm-m+1}\right\rfloor$.
It follows that either $\sigma(i-1,0)=((i-1)m,\overline{0})$ and
$\sigma(i,0)=(im,0)$
or $\sigma(i-1,0)=((i-1)m,0)$ and $\sigma(i,0)=(im,\overline{0})$.
In either case we see that $\sigma(\widetilde{\tau}(i,0))=\tau^m(\sigma(i,0))$.

\noindent (iii) Finally, suppose that $i=0$. Then $i-1\equiv (n-1)m\mod nm-m+1$, and,
as in Case (b)(i),
$\left\lfloor\frac{im}{nm-m+1}\right\rfloor=0$ is even and
$\left\lfloor\frac{(n-1)m}{nm-m+1}\right\rfloor=m-1$,
which means in this case that it is also even.
It follows that $\sigma(i-1,\overline{0})=((i-1)m,\overline{0})$ and
$\sigma(i,0)=(im,0)$ and thus that
$\sigma(\widetilde{\tau}(i,0))=\tau^m(\sigma(i,0))$.
\end{proof}

We therefore have:

\begin{theorem}\label{thm:D-component}
The translation quiver $\Gamma(D_n,m)$ can be realised as a connected
component of the restricted $m$th power of the translation quiver
$\Gamma(D_{nm-m+1},1)$.
\end{theorem}

Since $\mathcal{C}^m_{D_n}$ is equivalent to the additive hull of the
mesh category of $\Gamma(D_n,m)$, we obtain the following corollary.

\begin{cor}
The $m$-cluster category of type $D_n$ can be realised as the additive hull of
the mesh category of a connected component of the restricted $m$th power of
the Auslander-Reiten quiver of the cluster category of type $D_{nm-m+1}$.
For $m>2$ it is enough to take the usual $m$th power.
\end{cor}

\begin{figure}[ht]
\begin{center}
\includegraphics{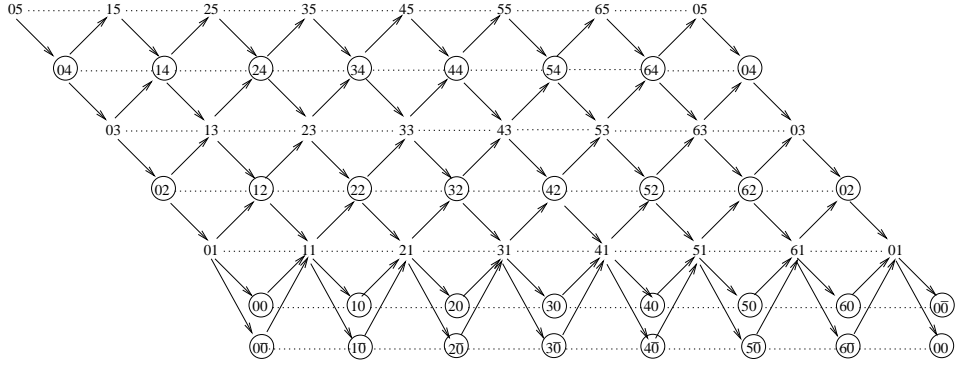} 
\caption{The translation quiver $\Gamma(D_7,1)$}
\label{fi:d7usual}
\end{center}
\end{figure}

\begin{figure}[ht]
\begin{center}
\includegraphics{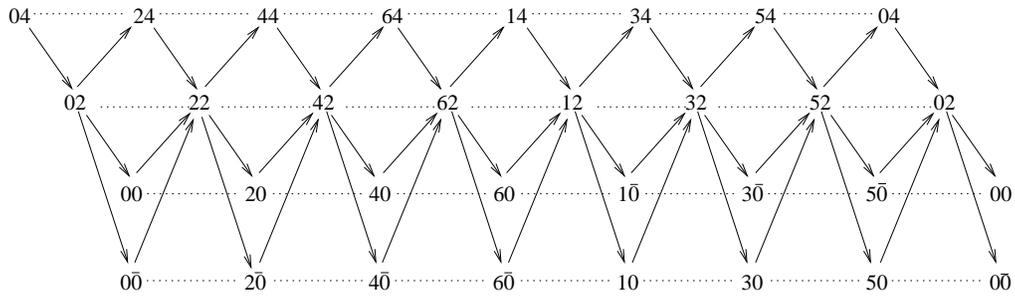} 
\caption{The connected component $\Gamma$ of the translation quiver
$\Gamma(D_7,1)$}
\label{fi:d4ind7}
\end{center}
\end{figure}

\begin{figure}[ht]
\begin{center}
\includegraphics{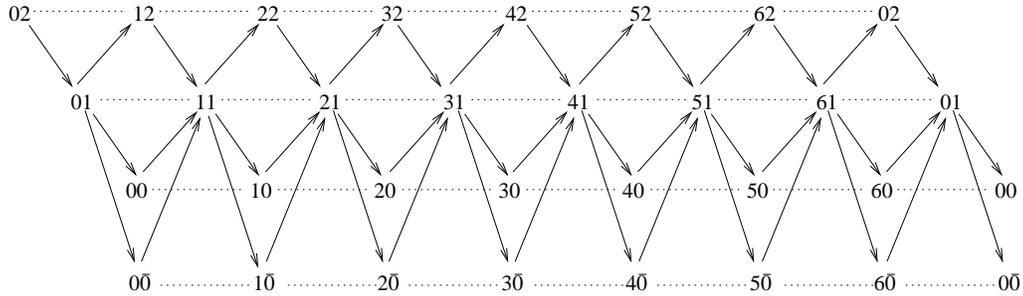} 
\caption{The translation quiver $\Gamma(D_4,2)$}
\label{fi:d4usual}
\end{center}
\end{figure}

\begin{example} \rm
We give an example of the theorem in the case where $n=4$ and $m=2$.
The theorem tells us that $\Gamma(D_4,2)$ is isomorphic to a connected
component of
$\mu_2(\Gamma(D_7,1))$. In Figure~\ref{fi:d7usual} we show the translation
quiver $\Gamma(D_7,1)$ with the vertices of $V=\im(\sigma')$ shown in circles.
In Figure~\ref{fi:d4ind7} we isolate the connected component
$\Gamma$ of $\mu_2(\Gamma(D_7,1))$ induced by $V$, and in
Figure~\ref{fi:d4usual} we indicate the
translation quiver $\Gamma(D_4,2)$ with the usual labelling of its vertices.
\end{example}

%
%
\section{Geometric realisation}

In this section, we give a geometric realisation of 
the $m$-cluster category of type $D_n$. To do so, we 
use certain $m$-arcs in a punctured $nm-m+1$-gon. Thus we 
are generalising the notion of tagged edges of 
Schiffler~\cite{schiffler} for the cluster 
category of type $D_n$ and the notion of 
$m$-diagonals of our 
work on $m$-cluster categories of type $A_n$~\cite{baurmarsh}.

Let $P$ be a punctured $N$-gon in the plane (later we shall specialise
to the case where $N=nm-m+1$). We label the vertices of $P$ clockwise. 
For $i\neq j\in\{1,2,\dots,N\}$, we denote by 
$B_{ij}$ the boundary path $i,i+1,\dots,j-1,j$, 
going clockwise around the boundary (taking the 
vertices $\mod N$). If $i=j$, we let $B_{ii}$ 
be the whole boundary path $i,i+1,\dots,i-1,i$ and 
$B_{ii}^{\bullet}$ denote the trivial path at $i$ 
consisting of the vertex $i$. 
The {\em length} $|B_{ij}|$ of the boundary path 
$B_{ij}$ is the number of vertices it runs through. 
Here, we count both the starting and end point unless 
$B_{ij}=B_{ii}^{\bullet}$. 
In particular, $|B_{ii}|=N+1$ and $|B_{ii}^{\bullet}|=1$.

\begin{figure}[ht] 
\begin{center}
\includegraphics{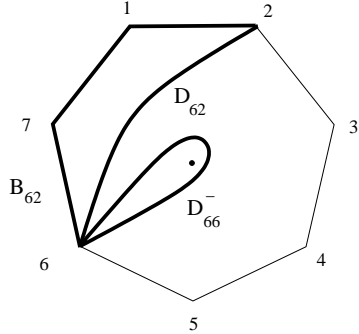} 
\caption{Punctured seven-gon with boundary path 
$B_{62}$ and arcs $D_{62}$ and $D_{66}^-$}\label{fig:7gon}
\end{center}
\end{figure}

As an example of a boundary path of length $4$, 
we have indicated $B_{62}$ 
inside a punctured $7$-gon in Figure~\ref{fig:7gon}. 

For $i\neq j$, and $j$ not the clockwise neighbour 
of $i$, an {\em arc} $D_{ij}$ is 
a line from $i$ to $j$ that is homotopic 
to the boundary path $B_{ij}$. 
If $j$ is the clockwise neighbour of $i$, there is 
no arc clockwise from $i$ to $j$ other than the boundary 
path $B_{ij}$. 
If $i=j$, we always tag the arc by $+$ or $-$, as 
in Schiffler's work~\cite{schiffler}. Such arcs are
denoted by $D_{ii}^+$ and $D_{ii}^-$. We will occasionally 
write $D_{ij}^{\pm}$ to denote an arbitrary arc and call 
it a {\em tagged arc}. In that 
case, if $i\neq j$, then $D_{ij}^{\pm}$ will only 
stand for the arc $D_{ij}$. 
As an example, the arcs $D_{62}$ and $D_{66}^-$ of a 
punctured $7$-gon are pictured in Figure~\ref{fig:7gon}. 
Addition of subscripts for the $D_{ij}$ will always be modulo $N$.

In what follows, we will use a slightly generalised 
version of a polygon. 
We will allow arcs $D_{ij}^{\pm}$ and sides $B_{i,i+1}$ of the polygon $P$
as sides of a polygon. We will say that such
a (generalised) polygon is {\em degenerate} if it 
has more sides than vertices. Note that such polygons may
or may not contain the puncture.

In the remainder, we will in particular be interested in 
the following types of generalised polygons and 
generalised degenerate polygons obtained from the 
regular $N$-gon $P$. 

{\bf Type (i):} A combination of an arc $D_{ij}$ with the boundary 
path $B_{ij}$, or of an arc $D_{ij}$ with the boundary 
path $B_{ji}$, $i\neq j$,
where in the former case, $2<|B_{ij}|\le N$ and in the 
latter case, $1<|B_{ji}|<N$. Such a polygon has 
$|B_{ij}|$ vertices (respectively, $|B_{ji}|$ vertices). 

{\bf Type (ii):} A combination of two arcs $D_{ij}$, $D_{ik}$ with 
the boundary path $B_{jk}$, or of $D_{ik}$, $D_{jk}$ 
with $B_{ij}$, where $i$, $j$, $k$ are all distinct 
and lie in clockwise order on $P$. 
Furthermore, in the former case, 
$1<|B_{jk}|<N-1$ and in the latter case, 
$1<|B_{ij}|<N-1$. Such a polygon has $|B_{jk}|+1$ vertices 
(respectively, $|B_{ij}|+1$ vertices). 

{\bf Type (iii):} A combination of an arc $D_{ii}^{\pm}$ 
with the boundary path $B_{ii}$, or with 
$B_{ii}^{\bullet}$. In the first case,  
the polygon is has $N+1$ sides and 
$N$ vertices. In the latter 
case, the polygon has one side and one vertex. 

{\bf Type (iv):} A combination of an arc $D_{ii}^{\pm}$ 
with an arc $D_{ij}$ and the boundary path 
$B_{ji}$, or 
a combination of $D_{ii}^{\pm}$ with a boundary path $B_{ij}$ 
and the arc $D_{ji}$ (where we always have 
$i\neq j$). In the former case, $1<|B_{ji}|<N$ and
in the latter case, $1<|B_{ij}|<N$.
Such a polygon has $|B_{ji}|+1$ sides and $|B_{ji}|$ 
vertices (respectively, $|B_{ij}|+1$ sides and $|B_{ij}|$ 
vertices). 

\begin{figure}[ht] 
\begin{center}
\includegraphics{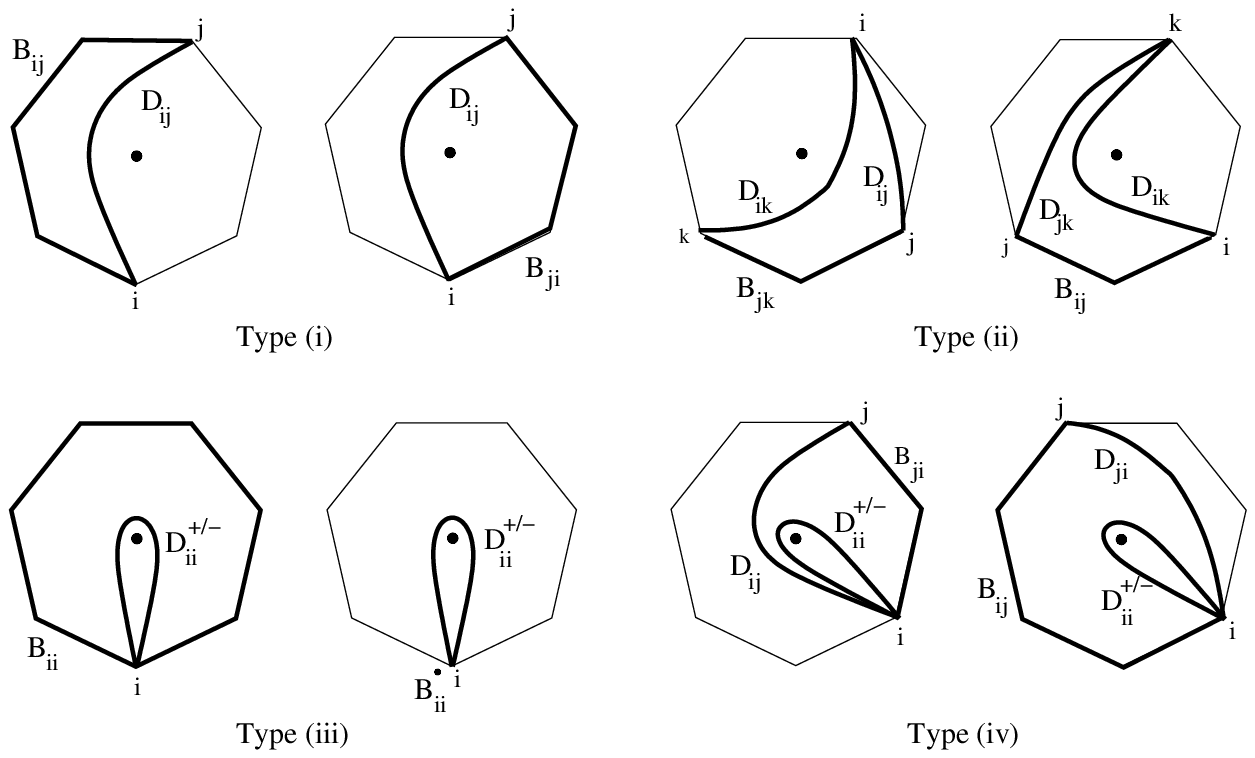} 
\caption{Generalised polygons (some degenerate)}\label{fig:polygons}
\end{center}
\end{figure}

Note that we can view type (iii) as the limit 
$j\mapsto i$ of type (i) 
and type (iv) as the limit $k\mapsto i$ or $j\mapsto k$ 
of type (ii). 
We show each of these four types in Figure~\ref{fig:polygons}. 

\begin{definition}\label{def:m-arc}
Let $D_{ij}^{\pm}$ be an arc of $P$. If $i\neq j$, we say 
that $D_{ij}$ is an {\em $m$-arc} if the following 
hold: 

(i) $D_{ij}$ and $B_{ij}$ form a $km+2$-gon for some $k$,

(ii) $D_{ij}$ and $B_{ji}$ form a $lm+1$-gon for 
some $l$. 

If $i=j$, $D_{ii}^{\pm}$ is a {\em tagged $m$-arc} 
if $D_{ii}^{\pm}$ and $B_{ii}$ form a degenerate 
$km+2$-gon for some $k$. 
\end{definition}
The parts (i) and (ii) in the definition of an $m$-arc 
also apply to  
the case $i=j$ if we use the 
boundary paths $B_{ii}$ and $B_{ii}^{\bullet}$. Namely, 
$D_{ii}^{\pm}$ is a tagged $m$-arc if $D_{ii}^{\pm}$ and $B_{ii}$ 
form a degenerate $km+2$-gon for some $k$ and if 
$D_{ii}^{\pm}$ and $B_{ii}^{\bullet}$ form a 
degenerate $1$-gon. 

\begin{example}
Let $P$ be a punctured $7$-gon and $m=2$. The arc 
$D_{62}$ is a $2$-arc (cf. Figure~\ref{fig:7gon}), since
the arc $D_{62}$ together 
with the boundary path $B_{62}$ forms a $4$-gon (i.e. $k=1$)
whereas $D_{62}$ and $B_{26}$ form a $5$-gon (i.e. $l=2$). 
Each of the arcs $D_{66}^{\pm}$ forms an $8$-gon together with 
$B_{66}$, and thus also is a $2$-arc.
\end{example}

We now define $m$-moves generalising the 
$m$-rotation for type $A_n$ of~\cite{baurmarsh} and 
the elementary moves for type $D_n$ of~\cite{schiffler}. 

\begin{definition}
Let $P$ be a punctured $N$-gon. An {\em $m$-move} arises
when there are two arcs in $P$ with a common end-point such that
the two arcs and a part of the boundary bound an unpunctured $m+2$-gon,
possibly degenerate. If the angle from the first arc to the
second at the common end-point is negative (i.e. clockwise), then we
say that there is an $m$-move taking the first arc to the second.
More precisely, it is a move of one of the following forms: 

(i) $D_{ij}\to D_{ik}$ if $D_{ij}$, $B_{jk}$ and $D_{ik}$ 
form an $m+2$-gon, $|B_{jk}|=m+1$. 

(ii) $D_{ij}\to D_{kj}$ if $D_{ij}$, $B_{ik}$ and $D_{kj}$ 
form an $m+2$-gon, $|B_{ik}|=m+1$. 

(iii) $D_{ij}\to D_{ii}^{\pm}$ if $D_{ij}$, $D_{ii}^{\pm}$ and $B_{ji}$ 
form a degenerate $m+2$-gon, $|B_{ji}|=m+1$. 

(iv) $D_{ii}^{\pm}\to D_{ji}$ if $D_{ii}^{\pm}$, $D_{ji}$ and $B_{ij}$ 
form a degenerate $m+2$-gon; $|B_{ij}|=m+1$. 
\end{definition}

\begin{figure}[ht] 
\begin{center}
\includegraphics{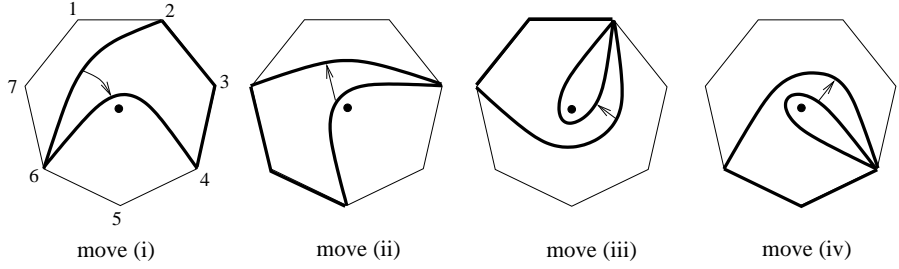} 
\caption{$2$-moves inside a heptagon}\label{fig:m-moves}
\end{center}
\end{figure}

In Figure~\ref{fig:m-moves}, we illustrate the four types 
of $m$-moves inside a heptagon, i.e. $n=4$, $m=2$. 

Our goal is to model the $m$-cluster category 
$\mathcal{C}^m(D_n)$ geometrically. 
To do so, we will from now on assume that $N=nm-m+1$, so the polygon $P$ 
has $nm-m+1$ vertices.

\begin{remark}
Let $P$ be a punctured polygon with $nm-m+1$ vertices and
let $i\neq j$. 
Then the two conditions of Definition~\ref{def:m-arc} 
are equivalent, i.e. $D_{ij}$ and $B_{ij}$ 
form a $km+2$-gon for some $k$ if and only if 
$D_{ij}$ and $B_{ji}$ form an $lm+1$-gon for some $l$. 
\end{remark}

We are now ready to define a translation quiver using  
the punctured polygon, $P$.

Let $\Gamma_{\odot}=\Gamma_{\odot}(n,m)$ be the quiver 
whose vertices are the tagged $m$-arcs of $P$ and whose 
arrows are given by $m$-moves. 
Let $\tau_m$ be the map sending an arc $D_{ij}^{\pm}$ 
to $D_{i-m,j-m}^{\pm}$ if $i\neq j$ or $m$ is even. 
If $i=j$ and $m$ is odd, we set 
$\tau_m(D_{ii}^{\pm})=D_{i-m,i-m}^{\mp}$. 
In other words, if $i\neq j$ or $m$ is even, 
$\tau_m$ rotates a tagged arc anti-clockwise around the center. 
In case $i=j$ and $m$ is odd, $\tau_m$ rotates the tagged 
arc anti-clockwise 
around the center and changes its tag. 

\begin{figure}[ht]
\begin{center}
\includegraphics{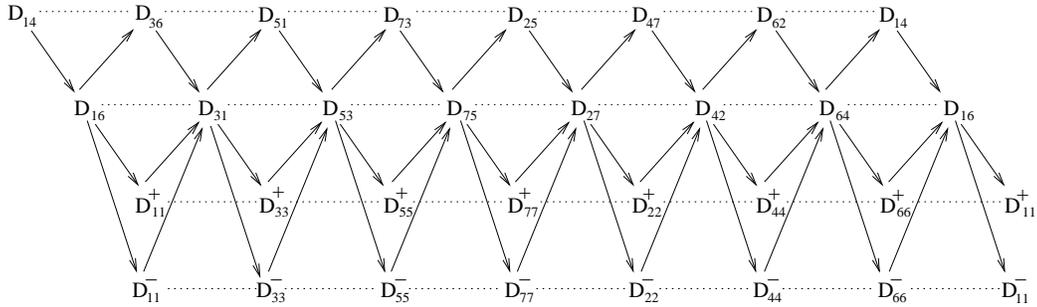} 
\caption{The translation quiver $\Gamma_{\odot}(4,2)$}
\label{fig:geom-4-2}
\end{center}
\end{figure}

Figure~\ref{fig:geom-4-2} shows the example $\Gamma_{\odot}(4,2)$ 
(we will see shortly, cf. Theorem~\ref{thm:odot-transl},
that $\Gamma_{\odot}(n,m)$ is a translation quiver).

\begin{lemma}\label{lm:sect-paths}
The sectional paths of length $m$ in 
$\Gamma(D_{nm-m+1},1)$ are of the 
form

\noindent
(i) $(i,j)\to (i,j-1)\to\dots\to (i,j-m)\ $ if $j-m> 0$, 

\noindent
(ii) 
$(i,j)\to (i,j-1)\to\dots\to (i,0)\ $ and $
(i,j)\to \dots\to (i,\overline{0})\ $ if $j=m$,

\noindent
(iii) $(i,j)\to\dots\to(i+m,j+m)\ $ if $j>0$ and 
$(i+m,j+m)$ exists,

\noindent
(iv) 
$(i,0)\to\dots\to (i+m,m)$ and 
$(i,\overline{0})\to\dots\to (i+m,m)\ $ 
if $(i+m,m)$ exists.
\end{lemma}
\begin{proof}
For (i) and (iii) we can 
argue as in~\cite[7.1]{baurmarsh}, using 
the vertices 
$\mathbb{Z}_{nm-m+1}\times\{0,\overline{0},1,\dots,n-2\}$ of 
$\Gamma(D_{nm-m+1},1)$ instead.  
The other cases follow with the same argument, using the 
assumption that the paths are restricted, 
i.e. excluding the sectional paths 
$(i,0)\to (i+1,1)\to (i+1,\overline{0})$ and 
$(i,\overline{0})\to (i+1,1)\to (i+1,0)$. 
\end{proof}

\begin{theorem} \label{thm:odot-transl}
The quiver $\Gamma_{\odot}$ is a translation quiver isomorphic to
to the Auslander-Reiten quiver of ${\mathcal C}^m_{D_n}$.
\end{theorem}

\begin{proof}
It is enough to show that $\Gamma_{\odot}$ is isomorphic to the image
$\Gamma$ of the map $\sigma$ from Section~\ref{se:mthpower}
and that, under the isomorphism, the map $\tau_m$ on
$\Gamma_{\odot}$ corresponds to $\tau^m$ on $\Gamma$.

Recall that the vertices of 
$\Gamma$ are $V:=\{(r,s)\in V(D_{N}):\ m|s\}$ 
(using the convention that $m$ divides $\overline{0}$), recalling
that $N=nm-m+1$.
In other words, $V$ is the subset in 
$V(D_{N})=\mathbb{Z}_{N}
\times\{0,\overline{0},1,\dots,N\}$ of 
the vertices whose second coordinate is divisible by $m$. 

We will now define a map $\rho:\ V(\Gamma_{\odot})\to\ V$,
where $V(\Gamma_{\odot})$ denotes the set of vertices of $\Gamma_{\odot}$.
Note that $m$-arcs in $V(\Gamma_{\odot})$ going through two
distinct vertices are always of form $D_{i,i+1+km}$. On such $m$-arcs,
$\rho$ is defined as follows: 
\[
\rho(D_{i,i+1+km})=(lm,(n-1-k)m)\in\mathbb{Z}_{N}
\times\{0,\overline{0},1,\dots,N-2\}
\] 
where $i\equiv lm+1$ 
modulo $N$ and $k=1,\dots,n-2$.  

\noindent
On arcs $D_{ii}^{\pm}$, $\rho$ is defined as follows. 
\begin{eqnarray*}
\rho(D_{ii}^{+}) & = & \left\{ \begin{array}{cc}
(lm,0), & \text{if $i$ is odd}, \\
(lm,\overline{0}), & \text{if $i$ is even}.\end{array}\right. \\
\\
\rho(D_{ii}^{-}) & = & \left\{ \begin{array}{cc}
(lm,\overline{0}), & \text{if $i$ is odd}, \\
(lm,0), & \text{if $i$ is even}.\end{array}\right. \\
\end{eqnarray*}
(where $i\equiv lm+1$ modulo $N$). 

To see that $\rho$ is a bijection, we divide $V(\Gamma_{\odot})$ 
up into $n$ types of arcs. Let 
$V_1$ be the set of arcs of the form $D_{i,i+m+1}$ ($i=1,\dots,N$), 
i.e. the arcs homotopic to a boundary path of length $m+2$. 
Then $\rho$ sends each element of $V_1$ to a vertex of the 
top row of $V$, 
$D_{i,i+m+1}\mapsto(lm,(n-2)m)$ (where $lm\equiv i-1\mod N$). 
It is straightforward to check that $\rho$ induces a bijection 
from the set $V_1$ to the top row of $V$.  

More generally, for $k=1,\dots,n-2$, let 
$V_k$ be the set of arcs of the form $D_{i,i+km+1}$, i.e. the 
set of arcs homotopic to a boundary path of length $km+2$. 
Since 
$\rho(D_{i,i+km+1})=(lm,(n-1-k)m)$, $\rho$ sends the arcs in 
$V_k$ to the $k$th row (from the top) of $V$ ($k\le n-2$). 
Clearly, this is also a bijection. 

Furthermore, the arcs $D_{ii}^{\pm}$ are sent to the two last rows 
of $V$, also bijectively. 
Thus, we have that $\rho$ is a bijection from 
$V(\Gamma_{\odot})$ to $V$. 

Next, we
observe that the arrows given by the $m$-moves
are the same as the arrows in $\Gamma$: 
for arcs in $V_k$ with $1\leq k<n-2$, an 
$m$-move sends $D_{i,i+1+km}$ to $D_{i,i+1+(k+1)m}$ 
or $D_{i,i+1+km}$ to $D_{i+m,i+1+km}$
whereas a restricted sectional path of length $m$ 
sends $(lm,(n-1-k)m)$ to $(lm,(n-2-k)m)$ (type (i) in 
Lemma~\ref{lm:sect-paths}) or to 
$((l+1)m,(n-k)m)$ (type (iii) in Lemma~\ref{lm:sect-paths}). 
For arcs in $V_{n-2}$, an $m$-move sends 
$D_{i,i+1+(n-1)m}$ to $D_{ii}^+$, to $D_{ii}^-$ or to 
$D_{i+m,i+1+(n-1)m}$ whereas 
a restricted sectional path of length $m$ sends $(lm,m)$ 
to $(lm,0)$, to $(lm,\overline{0})$ or to 
$((l+1)m,2m)$ (types (ii) and (iii) in Lemma~\ref{lm:sect-paths}).
Finally, arcs $D_{ii}^{\pm}$ are sent to $D_{i+m,i}$ by $m$-moves, 
and restricted sectional paths of length $m$ send 
$(lm,0)$ to $((l+1)m,m)$ and $(lm,\overline{0})$ to $((l+1)m,m)$
(type (iv) in Lemma~\ref{lm:sect-paths}).

Furthermore, the translation maps correspond: on 
$V_k$ (with $1\leq k\le n-2$)
$\tau_m(D_{i,i+1+km})=D_{i-m,i+1+(k-1)m}$ (subscripts taken 
$\mod N$) and on the $n-2$ first rows from the top, 
$\tau^m(lm,(n-1-k)m)=((l-1)m,(n-1-k)m)$ (first entries taken $\mod N$). 

If $i>1$ then $\tau_1(D_{ii}^{\pm})=D_{i-1,i-1}^{\mp}$ while
$\tau(i-1,0)=(i-2,0)$ and $\tau(i-1,\overline{0})=(i-2,\overline{0})$.
If $i=1$ then $\tau_1(D_{11}^+)=D_{N,N}^-$ while
$$\tau(0,0)=\left\{ \begin{array}{cc} (N-1,\overline{0}), & \text{if $N$ is odd}, \\
(N-1,0), & \text{if $N$ is even.} \\ \end{array}\right.$$
It follows that $\tau(\rho(D_{ii}^+))=\rho(\tau_1(D_{ii}^+))$
for all $i$. A similar argument applies to the tagged arcs $D_{ii}^-$.
Since $\tau_m=\tau_1^m$, we see that $\tau^m(\rho(D_{ii}^{\pm}))=
\rho(\tau_m(D_{ii}^{\pm}))$ for all $i$.
We have seen that $\rho$ induces an isomorphism of quivers and
$\tau^m\rho(D_{ij}^{\pm})=\rho(\tau_m(D_{ij}^{pm}))$ for all arcs
$D_{ij}^{\pm}$. It follows that $\Gamma_{\odot}$ is a translation quiver
and that $\rho$ is an isomorphism of translation quivers.
\end{proof}

%
%

\section{A toral translation quiver} \label{se:toralexample}

\begin{figure}[ht]
\begin{center}
\includegraphics{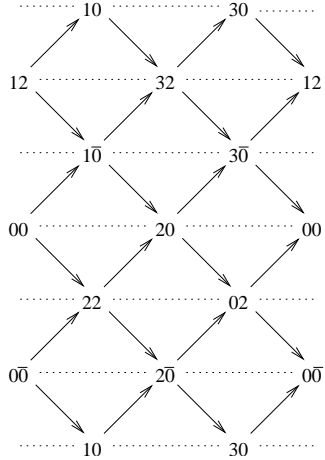} 
\caption{A connected component of the translation quiver 
$\Gamma(D_4,1)^2$}
\label{fi:toralquiver}
\end{center}
\end{figure}

In this section we give an example of a toral translation quiver arising
from the cluster category $\mathcal{C}_{D_4}^1$ of type $D_4$.
The Auslander-Reiten quiver of $\mathcal{C}_{D_4}^1$,
$\Gamma(D_4,1)$, is shown in
Figure~\ref{fi:quiverd4}. A connected component of its
(unrestricted) square, $\Gamma(D_4,1)^2$ is shown in
Figure~\ref{fi:toralquiver}. The underlying
topological space $|\Gamma(D_4)^2|$ (in the sense of Gabriel and
Riedtmann; see~\cite[p51]{ringel}) is a torus.


%
%

\section{The components of $\mu_m(\Gamma(D_n,1))$}

We have seen in Theorem~\ref{thm:D-component} that the Auslander-Reiten quiver
of the $m$-cluster category of type $D_n$ is a connected 
component of the restricted $m$-th power \linebreak
$\mu_m(\Gamma(D_{nm-m+1},1))$.
In this section, 
we describe the other components arising in the 
restricted $m$-th power 
of the translation quiver $(\Gamma(D_{nm-m+1},1),\tau)$. 

\begin{prop}\label{prop:A-components}
The quiver $\mu_m(\Gamma(D_{nm-m+1},1))$ has $m-1$ connected 
components isomorphic to the Auslander-Reiten quiver
of $D^b(A_{n-1})/\tau^{nm-m+1}$. 
\end{prop}
\begin{proof}
We consider the following subset of the vertices of the quiver 
\linebreak $\mu_m(\Gamma(D_{nm-m+1},1))$: 
\begin{eqnarray*}
X_k & := & \{(i,j)\mid i\in\mathbb{Z}_{nm-m+1}, j\equiv k\mod m\}\\
 & = & \mathbb{Z}_{nm-m+1}\times\{k,m+k,\dots,(n-2)m+k\}. 
\end{eqnarray*}
Such a set $X_k$ is a union of rows in the quiver
$\mu_m(\Gamma(D_{nm-m+1},1))$. 
We show that for each $1\le k\le m-1$, the translation quiver 
generated by $X_k$ (i.e. the full subquiver induced by $X_k$, 
together with $\tau^m$) is a connected 
component of $\mu_m(\Gamma(D_{nm-m+1},1))$. This is done by first showing
that $X_k$ is closed under $\tau^m$ and under taking 
restricted sectional paths of length $m$. This tells us that 
$X_k$ is a union of connected components of $\mu_m(\Gamma(D_{nm-m+1},1))$. 
Then we show that $X_k$ is connected, hence is a single 
component. \vspace{.2cm}

1) The set $X_k$ is closed under the translation $\tau^m$, since,
by definition, $X_k$ is the union of all vertices of certain 
rows and $\tau^m$ shifts vertices along a row.

2) The set $X_k$ is closed under restricted sectional paths of length 
$m$: we have seen that these paths are of the 
form $(ij)\to\dots\to(i,j-m)$ 
or $(ij)\to\dots\to (i+m,j-m)$, cf. Lemma~\ref{lm:sect-paths}. 
In particular, the new second entry is still congruent to $k$ modulo $m$. 

3) The subset $X_k$ is connected: note that $m$ is coprime to 
$nm-m+1$. Hence, the $\tau^m$-orbit of any vertex $(i,j)$ 
($i\in\mathbb{Z}_{nm-m+1}, j\equiv k\mod m$) is the same as 
the $\tau$-orbit of $(i,j)$. In other words, we 
can use $\tau^m$ to get everywhere in any given 
row of $X_k$, in particular in the row through $(0,k)$. 
Using the arrows and starting at $(0,k)$, 
we can get to any other row of $X_k$. 
\vspace{.2cm}
Now by definition, $X_k$ is the union of $n-1$ rows, namely the 
rows $(\cdot,k)$, $(\cdot,m+k)$ up to $(\cdot,(n-2)m+k)$. 
Each row is of length $nm-m+1$. 
It is clear from the arrows in $X_k$ that $X_k$ is isomorphic to
the Auslander-Reiten quiver of $D^b(A_{n-1})/\tau^{nm-m+1}$. 
\end{proof}

Thus we obtain a complete description of the 
restricted $m$-th power of \linebreak $\Gamma(D_{nm-m+1},1)$.

\begin{theorem}
The restricted $m$-th power $\mu_m(\Gamma(D_{nm-m+1},1),\tau^m)$ 
is the union of the following connected components: 
\[
\mu_m(\Gamma(D_{nm-m+1},1),\tau^m)=\Gamma_{\odot}(n,m)\cup
\bigcup_{k=i}^{m-1} \Gamma(D^b(A_{n-1})/\tau^{nm-m+1}),
\]
where $\Gamma(D^b(A_{n-1})/\tau^{nm-m+1})$ denotes the Auslander-Reiten
quiver of \linebreak $D_b(A_{n-1})/\tau^{nm-m+1}$.
\end{theorem}
\begin{proof}
The statement follows from 
Theorem~\ref{thm:D-component},
Proposition~\ref{prop:A-components} and the observation 
that the vertices of $(\mu_m(\Gamma(D_{nm-m+1},1)),\tau^m)$ 
are exhausted by the subsets $X_k$ ($k=1,\dots,m-1$) together 
with the vertices of $\Gamma_{\odot}(n,m)$. 
\end{proof}
%
%
%

\end{document}